\documentclass[leqno]{article}

\newtheorem{theorem}{Theorem}
\newtheorem{lemma}[theorem]{Lemma}
\newtheorem{proposition}[theorem]{Proposition}
\newtheorem{definition}[theorem]{Definition}
\newtheorem{corollary}[theorem]{Corollary}

\newcommand{\begintheorem}{\addtocounter{equation}{1}\begin{theorem}}
\newcommand{\beginlemma}{\addtocounter{equation}{1}\begin{lemma}}
\newcommand{\beginproposition}{\addtocounter{equation}{1}\begin{proposition}}
\newcommand{\begindefinition}{\addtocounter{equation}{1}\begin{definition}}
\newcommand{\begincorollary}{\addtocounter{equation}{1}\begin{corollary}}

\begin{document}

\title{Some topics in complex and harmonic analysis, 5}

\author{Stephen William Semmes	\\
	Rice University		\\
	Houston, Texas}

\date{}

\maketitle

	Let $f(x)$, $\phi(x)$ be continuous functions on ${\bf R}^n$,
and assume for simplicity that $f$ is bounded and that $\phi$ is
integrable.  For instance, $\phi$ might satisfy a bound like
\begin{equation}
	|\phi(x)| \le \frac{C}{1 + |x|^{n+1}}
\end{equation}
for some nonnegative real number $C$ and all $x \in {\bf R}^n$.  In
this event we can define the convolution of $f$ and $\phi$ in the
usual manner,
\begin{equation}
	(f * \phi)(x) = \int_{{\bf R}^n} f(y) \, \phi(x - y).
\end{equation}
This is especially simple if $\phi$ has compact support, in which
case $f$ can be any continuous function on ${\bf R}^n$.

	Let us make the normalizing assumption that the integral
of $\phi$ on ${\bf R}^n$ is equal to $1$.  For each positive real
number $t$ define $\phi_t(x)$ by
\begin{equation}
	\phi_t(x) = t^{-n} \, \phi(t^{-1} \, x).
\end{equation}
Thus $\phi_t$ is a continuous integrable function on ${\bf R}^n$ for
all $t > 0$ whose integral is also equal to $1$.  For small $t$
$\phi_t$ is basically mostly concentrated near $0$, while for $t$
large $\phi_t$ is more diffuse.  In particular, for each $r > 0$ the
integral of $\phi_t$ on the ball with center $0$ and radius $r$ tends
to $1$ as $t \to 0$.

	Using the continuity of $f$ one can check that $(\phi_t *
f)(x)$ tends to $f(x)$ as $t \to 0$ for all $x \in {\bf R}^n$.
Basically, $(f * \phi_t)(x)$ is an average of $f$ which is mostly
concentrated around $x$ as $t \to 0$, and thus it tends to $f(x)$ as
$t \to 0$.  Because $f$ is uniformly continuous on compact subsets of
${\bf R}^n$, one can show that $f * \phi_t$ converges to $f$ uniformly
on compact subsets of ${\bf R}^n$ as $t \to 0$.  If $f$ is uniformly
continuous, then $f * \phi_t$ converges to $f$ uniformly on ${\bf
R}^n$.  There are analogous statements for unbounded functions $f$
under suitable conditions on $\phi$.

	In some cases it may be that the limit of $f * \phi_t(x)$ as
$t \to 0$ exists even though $f$ is not continuous at $x$.  As a basic
scenario, suppose that $n = 1$, and that the right and left limits of
$f$ exist at a point $x$.  Suppose also that $\phi$ is an even
function, which is to say that $\phi(-x) = \phi(x)$.  In this case one
can check that $(f * \phi_t)(x)$ tends to the average of the left and
right limits of $f$ at $x$ as $t \to 0$.

	Let us continue to suppose that $n = 1$, and consider the case
where $\phi$ is a rational function on the real line.  In other words,
$\phi(x)$ can be written as $p(x)/ q(x)$, where $p$, $q$ are
polynomials and $q$ does not vanish on the real line.  In order
for $\phi$ to be integrable we should assume that the degree of $q$
is at least the degree of $p$ plus $2$.

	Notice that $\phi_t$ is then a rational function as well for
all $t > 0$.  Let us assume for simplicity that $f$ has compact
support in the real line.  The convolution $f * \phi_t$ is an integral
of translates of the rational function of $\phi$, and we can
approximate it by finite sums of translates of $\phi_t$, using Riemann
sums.  Of course finite sums of translates of a rational function are
again rational functions.

	Thus we get a nice way to approximate $f$ by rational
functions, namely by approximating $f$ by $f * \phi_t$ and then
approximating the convolution by a finite sum of translates of
$\phi_t$.  A rational function is in particular analytic,
which means that it has a convergent Taylor series expansion
in the neighborhood of any point.  The size of the nieghborhood
may be quite small, because the poles of the rational function
may be near by.

	Fix a closed and bounded interval $[a, b]$ in the real line,
and suppose that $r(x)$ is a rational function with no poles on this
interval.  Using partial fractions we can write $r(x)$ as a linear
combination of polynomials and rational functions of the form $(x +
c)^{-l}$, where $c$ is a complex number not in the interval $[a, b]$
and $l$ is a positive integer.  Each of these building blocks has a
convergent Taylor series expansion about some point in the complex
plane which converges uniformly on the interval $[a, b]$.  In this way
one can see that a rational function can be uniformly approximated
by polynomials on any closed and bounded interval in the real line
on which it does not have a pole.  One can use this and the earlier
arguments to show Weierstrass' approximation theorem, to the effect
that a continuous function on a closed and bounded interval in the
real line can be approximated uniformly by polynomials.

	In any dimension one can choose $\phi$ to be continuously
differentiable of all orders and to have compact support, so that if
$f$ is a continuous function on ${\bf R}^n$, then $f * \phi_t$ is
continuously differentiable of all orders and converges to $f$
unformly on compact subsets of ${\bf R}^n$, assuming that the integral
of $\phi$ is equal to $1$.  One might prefer to choose $\phi$ to be
real-analytic, so that $\phi$ has a convergent Taylor expansion in a
neighborhood of any point.  In this case the support of $\phi$ will be
all of ${\bf R}^n$, since otherwise $\phi$ would be identically equal
to $0$.  One can still choose $\phi$ to have enough decay to be
integrable, so that $f * \phi_t$ is defined for a suitable class of
functions $f$ and converges to $f$ uniformly on compact subsets of
${\bf R}^n$.

	Another interesting class of functions $\phi$ to consider
on ${\bf R}^n$ are functions of the form
\begin{equation}
	\phi(x) = \theta_1(x_1) \, \theta_2(x_2) \, \cdots \, \theta_n(x_n),
\end{equation}
where the $\theta_j$'s, $1 \le j \le n$, are continuous integrable
functions on the real line with integral equal to $1$.  Again $\phi_t$
would have the same form.  If the $\theta_j$'s have compact support
and $f$ is a continuous function on ${\bf R}^n$ with compact support,
then $f * \phi_t$ converges uniformly to $f$ on ${\bf R}^n$ and has
support contained in a fixed compact set when $t \le 1$, say, and $f *
\phi_t$ can be approximated by finite sums of products of functions of
one variable with supports contained in a fixed compact set, because
of the form of $\phi$.  In short we get a nice way to approximate $f$
uniformly by finite sums of products of functions of one variable with
supports contained in a fixed compact set.

	One might instead wish to choose functions $\phi$ which are
radial, which is to say that $\phi$ can be written as $\rho(|x|)$,
where $\rho$ is a continuous function of one variable.  This is the
same as saying that $\phi$ is invariant under orthogonal linear
transformations on ${\bf R}^n$, and of course $\phi_t$ is a radial
function for all $t > 0$ when $\phi$ is radial.  In this situation
the mapping which sends a function $f$ to $f * \phi_t$, under
suitable integrability conditions, has the nice feature that
it commutes with orthogonal linear transformations on ${\bf R}^n$.
In other words, if $R$ is an orthogonal linear transformation on
${\bf R}^n$, and one first sends $f$ to the composition $f \circ R$,
and then convolves the result with $\phi_t$, then that is the same
as taking $f * \phi_t$ and composing it with $R$.  In particular,
$f * \phi_t$ is a radial function if $f$ is radial.

\end{document}